\documentclass{amsart}

\theoremstyle{definition}

\theoremstyle{remark}

\reversemarginpar

\begin{document}

\title[Gamma function]{The integrals in Gradshteyn and Ryzhik. Part 4: \\
The gamma function.}

\author{Victor H. Moll}
\address{Department of Mathematics,
Tulane University, New Orleans, LA 70118}
\email{vhm@math.tulane.edu}

\subjclass{Primary 33}

\date{\today}


\begin{abstract}
We present a systematic derivation of some definite integrals in the classical 
table of Gradshteyn and Ryzhik that can be reduced to the gamma function.
\end{abstract}

\maketitle

\newcommand{\nn}{\nonumber}
\newcommand{\ba}{\begin{eqnarray}}
\newcommand{\ea}{\end{eqnarray}}
\newcommand{\ift}{\int_{0}^{\infty}}
\newcommand{\ione}{\int_{0}^{1}}
\newcommand{\ifft}{\int_{- \infty}^{\infty}}
\newcommand{\no}{\noindent}
\newcommand{\Ftwo}{{{_{2}F_{1}}}}
\newcommand{\realpart}{\mathop{\rm Re}\nolimits}
\newcommand{\imagpart}{\mathop{\rm Im}\nolimits}

\newtheorem{Definition}{\bf Definition}[section]
\newtheorem{Thm}[Definition]{\bf Theorem} 
\newtheorem{Example}[Definition]{\bf Example} 
\newtheorem{Lem}[Definition]{\bf Lemma} 
\newtheorem{Note}[Definition]{\bf Note} 
\newtheorem{Cor}[Definition]{\bf Corollary} 
\newtheorem{Prop}[Definition]{\bf Proposition} 
\newtheorem{Problem}[Definition]{\bf Problem} 
\numberwithin{equation}{section}

\section{Introduction} \label{intro} 
\setcounter{equation}{0}

The table of integrals \cite{gr} contains some evaluations that can be 
derived by elementary means from the {\em gamma function}, defined by 
\begin{equation}
\Gamma(a) = \ift x^{a-1} e^{-x} \, dx.
\label{gamma-def}
\end{equation}
\noindent
The convergence of the integral in (\ref{gamma-def}) requires $a > 0$. 
The goal of this paper is to present some of these evaluations 
in a systematic manner. The  
techniques developed here will be employed 
in future publications. The reader will find in \cite{irrbook} analytic 
information about this important function.

The gamma function represents the extension of factorials to real parameters.
The value 
\begin{equation}
\Gamma(n) = (n-1)!, \text{ for } n \in \mathbb{N} 
\end{equation}
\noindent
is elementary. On the other hand, the special value 
\begin{equation}
\Gamma \left( \tfrac{1}{2} \right) = \sqrt{\pi}
\end{equation}
\noindent
is equivalent to the well-known {\em normal integral}
\begin{equation}
\ift  \text{exp}(-t^{2}) \, dt =
 \tfrac{1}{2} \Gamma \left( \tfrac{1}{2} \right).
\label{normal}
\end{equation}
\noindent
The reader will find in \cite{irrbook} proofs of Legendre's duplication 
formula
\begin{equation}
\Gamma \left( x + \tfrac{1}{2} \right) = \frac{\Gamma(2x) \sqrt{\pi}}{\Gamma(x) 
\, 2^{2x-1}},
\label{legendre}
\end{equation}
\noindent
that produces for $x = m \in \mathbb{N}$ the values 

\begin{equation}
\Gamma \left( m + \tfrac{1}{2} \right) = \frac{\sqrt{\pi}}{2^{2m}} \, 
\frac{(2m)!}{m!}. 
\label{gammahalf}
\end{equation}
\noindent
This appears as $\mathbf{3.371}$ in \cite{gr}. 

\section{The introduction of a parameter} \label{sec-parameter} 
\setcounter{equation}{0}

The presence of a parameter in a definite integral provides great amount of 
flexibility. The 
change of variables $ x = \mu t$ in (\ref{gamma-def}) yields
\begin{equation}
\Gamma(a) = \mu^{a} \ift t^{a-1} e^{-\mu t} \, dt. 
\label{gamma1}
\end{equation}

\noindent
This appears as $\mathbf{3.381.4}$ in \cite{gr} and the choice $a = n+1$, with 
$n \in \mathbb{N}$, that reads
\begin{equation}
\ift t^{n}e^{- \mu t} \, dt = n! \, \mu^{-n-1}
\end{equation}
\noindent
appears as $\mathbf{3.351.3}$. \\

 The special case 
$a = m+ \tfrac{1}{2}$, that appears as $\mathbf{3.371}$ in \cite{gr}, 
yields
\begin{equation}
\ift t^{m - \tfrac{1}{2}} e^{- \mu t} \, dt = 
\frac{\sqrt{\pi}}{2^{2m}} \, 
\frac{(2m)!}{m!} \mu^{-m - \tfrac{1}{2}}, 
\end{equation}
\noindent
is consistent with (\ref{gammahalf}). 

\medskip

The combination 
\begin{equation}
\ift \frac{e^{- \nu x} - e^{- \mu x}}{x^{ \rho + 1} } \, dx = 
\frac{\mu^{\rho} - \nu^{\rho}}{\rho} \, \Gamma(1 - \rho),
\label{int-11}
\end{equation}
\noindent
that appears as $\mathbf{3.434.1}$ in \cite{gr} can 
now be evaluated directly. The parameters are restricted by convergence: 
$\mu, \, \nu > 0$ and $\rho < 1$. The 
integral $\mathbf{3.434.2}$
\begin{equation}
\ift \frac{e^{- \mu x} - e^{- \nu x}}{x} \, dx  = \ln \frac{\nu}{\mu},
\end{equation}
\noindent
is obtained from (\ref{int-11}) by passing to the limit as $\rho \to 0$. This 
is an example of {\em Frullani integrals} that will be discussed in a future 
publication. \\

The reader will be able to check $\mathbf{3.478.1}$:

\begin{equation}
\ift x^{\nu-1} \text{exp}(- \mu x^{p}) \, dx = \frac{1}{p} \mu^{-\nu/p} 
\Gamma \left( \frac{\nu}{p} \right), 
\end{equation}
\noindent
and $\mathbf{3.478.2}$:
\begin{equation}
\ift x^{\nu-1} \left[ 1 - \text{exp}(-\mu x^{p} ) \right] \, dx = 
-\frac{1}{|p|} \mu^{-\nu/p} \Gamma \left( \frac{\nu}{p} \right)
\end{equation}
\noindent
by introducing appropriate parameter reduction. \\

The parameters can be used to prove many of the classical identities for 
$\Gamma(a)$.

\begin{Prop}
The gamma function satisfies 
\begin{equation}
\Gamma(a+1) = a \, \Gamma(a). 
\end{equation}
\end{Prop}
\begin{proof}
Differentiate (\ref{gamma1}) with respect to $\mu$ to produce
\begin{equation}
0 = a \mu^{a-1} \ift t^{a-1}e^{-\mu t} \, dt - 
\mu^{a} \ift t^{a}e^{-\mu t} \, dt. 
\end{equation}
\noindent
Now put $\mu=1$ to obtain the result.
\end{proof}

\medskip

Differentiating (\ref{gamma-def}) with respect to the parameter $a$ yields
\begin{equation}
\Gamma'(a) = \ift x^{a-1} e^{-x} \, \ln x \, dx. 
\end{equation}
\noindent
Further differentiation introduces higher powers of $\ln x$: 
\begin{equation}
\Gamma^{(n)}(a) = \ift x^{a-1} e^{-x} \, \left( \ln x \right)^{n} \, dx. 
\end{equation}
\noindent
In particular, for $a=1$, we obtain:
\begin{equation}
\ift \left( \ln x \right)^{n} e^{-x} \, dx = \Gamma^{(n)}(1). 
\end{equation}
\noindent
The special case $n=1$ yields 
\begin{equation}
\ift e^{-x} \, \ln x \, dx = \Gamma'(1). 
\end{equation}
\noindent
The reader will find in \cite{irrbook}, page $176$ an elementary proof that 
$\Gamma'(1) = -\gamma$, where 
\begin{equation}
\gamma := \lim \limits_{n \to \infty} \sum_{k=1}^{n} \frac{1}{k} - \ln n 
\end{equation}
\noindent
is Euler's constant. This is one of the fundamental numbers of Analysis. 

On the other hand, differentiating (\ref{gamma1}) produces
\begin{equation}
\ift x^{a-1}e^{-\mu x} \left( \ln x \right)^{n} \, dx = 
\left( \frac{\partial}{\partial a} \right)^{n} \left[ \mu^{-a} \Gamma(a) 
\right], 
\label{derrr}
\end{equation}
\noindent
that appears as $\mathbf{4.358.5}$ 
in \cite{gr}. Using Leibnitz's differentiation 
formula we obtain
\begin{equation}
\ift x^{a-1}e^{-\mu x} \left( \ln x \right)^{n} \, dx = 
\mu^{-a} \sum_{k=0}^{n} (-1)^{k} \binom{n}{k} \left( \ln \mu \right)^{k} 
\Gamma^{(n-k)}(a).
\end{equation}
\noindent
In the special case $a=1$ we obtain
\begin{equation}
\ift e^{-\mu x} \left( \ln x \right)^{n} \, dx = 
\frac{1}{\mu} \sum_{k=0}^{n} (-1)^{k} \binom{n}{k} \left( \ln \mu \right)^{k} 
\Gamma^{(n-k)}(1).
\end{equation}
\noindent
The cases $n=1, \, 2, \, 3$ appear as $\mathbf{4.331.1}, \, 
\mathbf{4.335.1}$ and 
$\mathbf{4.335.3}$ respectively. 

In order to obtain analytic expressions for the terms $\Gamma^{(n)}(1)$, it
is convenient to introduce the {\em polygamma function} 
\begin{equation}
\psi(x) = \frac{d}{dx} \ln \Gamma(x).
\label{poly-def}
\end{equation}
\noindent
The derivatives of $\psi$ satisfy
\begin{equation}
\psi^{(n)}(x) = (-1)^{n+1} n! \, \zeta(n+1,x),
\label{der-psi}
\end{equation}
\noindent
where 
\begin{equation}
\zeta(z,q) = \sum_{n=0}^{\infty} \frac{1}{(n+q)^{z}}
\end{equation}
\noindent
is the {\em Hurwitz zeta function}. In particular this gives
\begin{equation}
\psi^{(n)}(1) = (-1)^{n+1} \, n! \, \zeta(n+1). 
\end{equation}
\noindent
The values of $\Gamma^{(n)}(1)$ can now be computed by recurrence via
\begin{equation}
\Gamma^{(n+1)}(1) = \sum_{k=0}^{n} \binom{n}{k} \Gamma^{(k)}(1) 
\psi^{(n-k)}(1),
\end{equation}
\noindent
obtained by differentiating $\Gamma'(x) = \psi(x) \Gamma(x)$.  \\

Using (\ref{der-psi}) the reader will be able to check the first few 
cases of (\ref{derrr}), we employ the notation $\delta = \psi(a) - \ln \mu$:
\noindent
\begin{eqnarray}
\ift x^{a-1} e^{-\mu x} \ln^{2}x \, dx & = & \frac{\Gamma(a)}{\mu^{a}}
\left\{ \delta^{2} + \zeta(2,a) \right\}, \nonumber  \\
\ift x^{a-1} e^{-\mu x} \ln^{3}x \, dx & = & \frac{\Gamma(a)}{\mu^{a}}
\left\{ \delta^{3} + 3 \zeta(2,a) \delta - 2 \zeta(3,a) \right\}, \nonumber \\
\ift x^{a-1} e^{-\mu x} \ln^{4}x \, dx & = & \frac{\Gamma(a)}{\mu^{a}}
\left\{ \delta^{4} + 6 \zeta(2,a) \delta^{2} - 8 \zeta(3,a) \delta +
3 \zeta^{2}(2,a) + 6 \zeta(4,a) \right\}. \nonumber 
\end{eqnarray}
\noindent
These appear as $\mathbf{4.358.2}, \, \mathbf{4.358.3}$ 
and $\mathbf{4.358.4}$, respectively.

\section{Elementary changes of variables} \label{sec-changes} 
\setcounter{equation}{0}

The use of appropriate changes of variables yields, from the basic 
definition (\ref{gamma-def}), the evaluation of 
more complicated definite integrals. For example, let $x = t^{b}$ to 
obtain, with $c = ab-1$,
\begin{equation}
\ift t^{c} \text{exp}(-t^{b}) \, dt = \frac{1}{b} \Gamma \left( \frac{c+1}{b} 
\right). 
\end{equation}
\noindent
The special case $a = 1/b$, that is $c=0$,  is
\begin{equation}
\ift  \text{exp}(-t^{b}) \, dt = \frac{1}{b} \Gamma \left( \frac{1}{b} \right),
\end{equation}
\noindent
that appears as $\mathbf{3.326.1}$ in \cite{gr}. The special case $b=2$ is the 
normal integral (\ref{normal}).

We can now introduce an extra 
parameter via $t = s^{1/b}x$. This produces 
\begin{equation}
\ift  x^{m} \text{exp}(- s x^{b}) \, dx = 
\frac{\Gamma(a)}{s^{a}b},
\label{gamma-4}
\end{equation}

\noindent
with $m = ab-1$. This formula appears (at least) three times in 
\cite{gr}:  $\mathbf{3.326.2}, \mathbf{3.462.9}$ 
and $\mathbf{3.478.1}$. Moreover, the case $s=1, \, 
c = (m + 1/2)n-1$ and $b=n$ appears as $\mathbf{3.473}$:
\begin{equation}
\ift \text{exp}(-x^{n}) x^{ \left(m+ \tfrac{1}{2} \right)n-1} \, dx = 
\frac{(2m-1)!!}{2^{m} \, n} \sqrt{\pi}. 
\end{equation}
\noindent
The form given here can be established using (\ref{gammahalf}). \\

Differentiating (\ref{gamma-4}) with respect to the parameter $m$ (keeping in
mind that $a = (m+1)/b$), yields

\begin{equation}
\ift x^{m} e^{-sx^{b}} \, \ln x \, dx = 
\frac{\Gamma(a)}{b^{2} \, s^{a}} \left[ \psi(a) - \ln s \right].
\label{gamma-7}
\end{equation}
\noindent
In particular, if $b=1$ we obtain
\begin{equation}
\ift x^{m} e^{-sx} \, \ln x \, dx = 
\frac{\Gamma(m+1)}{ s^{m+1}}  \left[ \psi(m+1) - \ln s \right].
\label{gamma-7a}
\end{equation}
\noindent
The case $m=0$ and $b=2$ gives
\begin{equation}
\ift e^{-sx^{2}} \, \ln x \, dx  = - \frac{\sqrt{\pi}}{4 \sqrt{s}} 
\left( \gamma + \ln 4s \right), 
\end{equation}
\noindent 
where we have used $\psi(1/2) = -\gamma - 2 \ln 2$. This appears 
as $\mathbf{4.333}$ 
in \cite{gr}. 

\medskip

An interesting example is $b=m=2$. Using the values 
\begin{equation}
\Gamma \left( \tfrac{3}{2} \right) = 
\sqrt{\pi}/2 \text{ and } \psi \left(\tfrac{3}{2} \right) = 
2 - 2 \ln 2 - \gamma 
\end{equation}
\noindent
the expression (\ref{gamma-7}) yields
\begin{equation}
\ift x^{2} e^{-sx^{2}} \, \ln x \, dx = 
\frac{1}{8s}(2 - \ln 4s - \gamma) \, \sqrt{\frac{\pi}{s}}. 
\label{form-99}
\end{equation}
\noindent
The values of $\psi$ at half-integers follow directly from (\ref{legendre}). 
Formula (\ref{form-99}) appears as $\mathbf{4.355.1}$ in 
\cite{gr}. Using (\ref{gamma-7}) it is easy to 
verify
\begin{equation}
\ift (\mu x^{2}- n)x^{2n-1}e^{-\mu x^{2}} \, \ln x \, dx = 
\frac{(n-1)!}{4 \mu^{n}}, 
\end{equation}
\noindent
and 
\begin{equation}
\ift (2 \mu x^{2} - 2n-1)x^{2n} e^{-\mu x^{2}} \, \ln x \, dx = 
\frac{(2n-1)!!}{2(2 \mu)^{n}} \sqrt{\frac{\pi}{\mu}}, 
\end{equation}
\noindent
for $n \in \mathbb{N}$. These
appear as, respectively,  $\mathbf{4.355.3}$ and 
$\mathbf{4.355.4}$ in \cite{gr}. The 
term $(2n-1)!!$ is the semi-factorial defined by
\begin{equation}
(2n-1)!! = (2n-1)(2n-3) \cdots 5 \cdot 3 \cdot 1.
\end{equation}

Finally, formula $\mathbf{4.369.1}$ in \cite{gr}
\begin{equation}
\ift x^{a-1} e^{-\mu x} \left[ \psi(a) - \ln x \right] \, dx = 
\frac{\Gamma(a) \, \ln \mu}{\mu^{a}}
\end{equation}
\noindent
can be established by the methods developed here.  The more ambitious 
reader will check that
\begin{equation}
\ift x^{n-1} e^{-\mu x} \left\{ [ \ln x - \tfrac{1}{2} \psi(n)]^{2} - 
\tfrac{1}{2} \psi'(n) \right\} \, dx = 
\frac{(n-1)!}{\mu^{n}} \left\{ [ \ln \mu - \tfrac{1}{2} \psi(n)]^{2} +
\tfrac{1}{2} \psi'(n) \right\}, \nonumber 
\end{equation}
\noindent
that is $\mathbf{4.369.2}$ in \cite{gr}.  \\

We can also write (\ref{gamma-7}) in the exponential scale to obtain
\begin{equation}
\int_{-\infty}^{\infty} t e^{mt} \text{exp}\left( -s e^{bt} \right) \, dt =
\frac{\Gamma(m/b)}{b^{2}s^{m/b}} \left( \psi \left( \frac{m}{b} \right)- 
\ln s \right). 
\end{equation}
\noindent
The special case $b=m=1$ produces
\begin{equation}
\int_{-\infty}^{\infty} te^{t} \, \text{exp}\left(-se^{t} \right) \, dt = 
- \frac{(\gamma + \ln s)}{s} 
\end{equation}
\noindent
that appears as $\mathbf{3.481.1}$. The second 
special case, appearing as $\mathbf{3.481.2}$, 
is $b=2,\, m=1$, that yields 
\begin{equation}
\int_{-\infty}^{\infty} te^{t} \, \text{exp}\left(-se^{2t} \right) \, dt = 
- \frac{\sqrt{\pi} \, (\gamma + \ln 4s)}{4 \sqrt{s}}. 
\end{equation}
\noindent
This uses the value $\psi(1/2) = -(\gamma + 2 \ln 2)$. \\

\medskip

There are many other possible changes of variables that lead to 
interesting evaluations. We conclude this section with one more:
let $x = e^{t}$ to convert (\ref{gamma-def}) into 
\begin{equation}
\int_{-\infty}^{\infty} \text{exp} \left(-e^{x} \right) \, e^{ax} \, dx = 
\Gamma(a). 
\end{equation}
\noindent
This is $\mathbf{3.328}$ in \cite{gr}. \\

As usual one should not prejudge the difficulty of a problem: the example 
$\mathbf{3.471.3}$ states that
\begin{equation}
\int_{0}^{a} x^{-\mu-1} (a-x)^{\mu-1} e^{-\beta/x} \, dx = 
\beta^{-\mu} a^{\mu-1} \Gamma(\mu) \text{ exp} \left(- \tfrac{\beta}{a} \right).
\end{equation}
\noindent
This can be reduced to the basic formula for the gamma function. Indeed, the 
change of  variables $t = \beta/x$ produces 
\begin{equation}
I = \beta^{-\mu}a^{\mu-1} \int_{\beta/a}^{\infty} 
\left( t - \beta/a \right)^{\mu-1} e^{-t} \, dt. 
\end{equation}
\noindent
Now let $y = t - \beta/a$ to complete the evaluation. The table \cite{gr} 
writes $\mu$ instead of $a$: it seems to be a bad idea to have $\mu$ and $u$ 
in the same formula, it leads to typographical errors that should be 
avoided. \\

Another simple change of variables gives the evaluation 
of $\mathbf{3.324.2}$:
\begin{equation}
\int_{-\infty}^{\infty} e^{-(x-b/x)^{2n}} \, dx = \frac{1}{n} 
\Gamma \left( \frac{1}{2n} \right). 
\end{equation}
\noindent
The symmetry yields
\begin{equation}
I = 2 \ift e^{-(x-b/x)^{2n}} \, dx. 
\end{equation}
\noindent
The change of variables $t=b/x$ yields, using $b>0$,
\begin{equation}
I = 2b \ift e^{-(t-b/t)^{2n}} \,  \frac{dt}{t^{2}}.
\end{equation}
\noindent
The average of these forms produces
\begin{equation}
I =  \ift e^{-(x-b/x)^{2n}} \, \left( 1 + \frac{b}{x^{2}} \right) \, dx.
\end{equation}
\noindent
Finally, the change of variables $u = x-b/x$ gives the result. Indeed, 
let $u = x-b/x$ and observe that $u$ is increasing when $b >0$. 
This restriction is missing in the table. Then we get
\begin{equation}
I = 2 \ift e^{-u^{2n}} \, du. 
\end{equation}
\noindent
This can now be evaluated via $v = u^{2n}$. \\

\noindent
{\bf Note}. In the case $b<0$ the change of variables $u=x-b/x$ has an 
inverse with two branches, splitting at $x= \sqrt{-b}$. Then we write
\begin{eqnarray}
I & := & 2 \ift e^{-(x-b/x)^{2n}} \, dx  \label{int11} \\
 & = & 2 \int_{0}^{\sqrt{-b}} e^{-(x-b/x)^{2n}} \, dx + 
       2 \int_{\sqrt{-b}}^{\infty} e^{-(x-b/x)^{2n}} \, dx. \nonumber
\end{eqnarray}
\noindent
The change of variables $u=x-b/x$ is now used in each of the integrals to 
produce
\begin{equation}
I = 2 \int_{2 \sqrt{-b}}^{\infty} 
\frac{u \, \text{exp}(-u^{2n}) \,du}{\sqrt{u^{2}+4b}}. 
\end{equation}
\noindent
The change of variables $z = \sqrt{u^{2}+4b}$ yields
\begin{equation}
I = 2 \ift \text{exp} \left( -(z^{2}-4b)^{n} \right). 
\end{equation}
\noindent
We are unable to simplify it any further.

\section{The logarithmic scale} \label{sec-logar} 
\setcounter{equation}{0}

Euler prefered the version 
\begin{equation}
\Gamma(a) = \int_{0}^{1} \left( \ln \frac{1}{u} \right)^{a-1} \, du. 
\label{euler-def}
\end{equation}
\noindent
We will write this as 
\begin{equation}
\Gamma(a) = \int_{0}^{1} \left( - \ln u \right)^{a-1} \, du,
\label{euler-def1}
\end{equation}
\noindent
for better spacing. 
Many of the evaluations in \cite{gr} follow this form. Section 
$\mathbf{4.215}$ in 
\cite{gr} consists of  four examples: 
the first one, $\mathbf{4.215.1}$ is 
(\ref{euler-def}) itself. The second one, labeled $\mathbf{4.215.2}$ and 
written as
\begin{equation}
\int_{0}^{1} \frac{dx}{\left( - \ln x \right)^{\mu}} = 
\frac{\pi}{\Gamma(\mu)} \text{cosec }\mu \pi, 
\end{equation}
\noindent
is evaluated as $\Gamma(1 - \mu)$ by (\ref{euler-def}).  The identity 
\begin{equation}
\Gamma(\mu) \Gamma(1 - \mu) = \frac{\pi}{\sin \pi \mu} 
\label{gamma-99}
\end{equation}
\noindent 
yields the given form. The reader will find in \cite{irrbook} a proof of 
this identity. The section concludes with the special values 
\begin{equation}
\int_{0}^{1} \sqrt{- \ln x} \, dx = \frac{\sqrt{\pi}}{2}, 
\end{equation}
\noindent
as $\mathbf{4.215.3}$ and  $\mathbf{4.215.4}$:
\begin{equation}
\int_{0}^{1} \frac{dx}{\sqrt{ - \ln x}} = \sqrt{\pi}. 
\end{equation}
\noindent
Both of them are special cases of (\ref{euler-def}). \\

The reader should check the evaluations $\mathbf{4.269.3}$:
\begin{equation}
\int_{0}^{1} x^{p-1} \, \sqrt{ -  \ln x } \, dx = \frac{1}{2}
\sqrt{\frac{\pi}{p^{3}}}, 
\end{equation}
\noindent
and $\mathbf{4.269.4}$:
\begin{equation}
\int_{0}^{1} \frac{x^{p-1} \, dx}{\sqrt{- \ln x }} = 
\sqrt{ \frac{\pi}{p}}
\end{equation}
\noindent
by reducing them to (\ref{gamma1}). Also $\mathbf{4.272.5}, \, 
\mathbf{4.272.6}$
and $\mathbf{4.272.7}$ 
\begin{eqnarray}
\int_{1}^{\infty} \left( \ln x \right)^{p} \frac{dx}{x^{2}} & = & 
\Gamma(1+p),  \\
\int_{0}^{1} \left( - \ln x \right)^{\mu -1} \, x^{\nu-1} \, dx 
& = & \frac{1}{\nu^{\mu}} \Gamma(\mu), \nonumber \\
\int_{0}^{1} \left( - \ln x \right)^{n - \tfrac{1}{2} } \, x^{\nu-1} 
\, dx & = & \frac{(2n-1)!!}{(2 \nu)^{n}} \sqrt{\frac{\pi}{\nu}},
\nonumber 
\end{eqnarray}
\noindent
can be evaluated directly in terms of 
the gamma function.  \\

Differentiating (\ref{euler-def}) with respect to $a$ yields 
$\mathbf{4.229.4}$ in \cite{gr}:
\begin{equation}
\int_{0}^{1} \ln \left( - \ln x \right) \, 
\left( - \ln x \right)^{a-1} \, dx = \Gamma'(a) = \psi(a) \Gamma(a), 
\end{equation}
\noindent
with $\psi(a)$ defined in (\ref{poly-def}). The special case $a=1$ is 
$\mathbf{4.229.1}$:
\begin{equation}
\int_{0}^{1} \ln \left( - \ln x \right) \, dx = - \gamma, 
\end{equation}
\noindent
and 
\begin{equation}
\int_{0}^{1} \ln \left( - \ln x \right) \, 
\frac{dx}{\sqrt{- \ln x }}  = -( \gamma + 2 \ln 2 ) \sqrt{\pi}, 
\end{equation}
\noindent 
that appears as $\mathbf{4.229.3}$, is obtained 
by using the values $\Gamma \left( \tfrac{1}{2} \right) 
= \sqrt{\pi}$ 
and $\psi \left( \tfrac{1}{2} \right) =  -(\gamma +  2 \ln 2)$.  \\

The same type of arguments confirms $\mathbf{4.325.11}$
\begin{equation}
\int_{0}^{1} \ln( - \ln x ) \, 
\frac{x^{\mu -1} \, dx}{\sqrt{- \ln x}} = 
-( \gamma + \ln 4 \mu) \sqrt{\frac{\pi}{\mu}}, 
\end{equation}
\noindent
and $\mathbf{4.325.12}$:
\begin{equation}
\int_{0}^{1} \ln ( - \ln  x ) \,  \left( - \ln x 
\right)^{\mu-1} \, x^{\nu-1} \, dx 
= \frac{1}{\nu^{\mu}} \Gamma(\mu) \left[ \psi(\mu) - \ln \nu \right]. 
\end{equation}
\noindent
In particular, when $\mu=1$ we obtain $\mathbf{4.325.8}$:
\begin{equation}
\int_{0}^{1} \ln ( - \ln x ) \,  x^{\nu-1} \, dx
= -\frac{1}{\nu} \left( \gamma + \ln \nu \right). 
\end{equation}

\section{The presence of fake parameters} \label{sec-fake} 
\setcounter{equation}{0}

There are many formulas in \cite{gr} that contain parameters. For example, 
$\mathbf{3.461.2}$ states that
\begin{equation}
\ift x^{2n} e^{-px^{2}} \, dx = \frac{(2n-1)!!}{2(2p)^{n}} \sqrt{\frac{\pi}{p}}
\end{equation}
\noindent
and $\mathbf{3.461.3}$ states that
\begin{equation}
\ift x^{2n+1} e^{-px^{2}} \, dx = \frac{n!}{2p^{n+1}}. 
\end{equation}
\noindent
The change of  variables $t = px^{2}$ eliminates the {\em fake} parameter $p$ 
and reduces $\mathbf{3.461.2}$ to 
\begin{equation}
\ift t^{n - \tfrac{1}{2} } e^{-t} \, dt = \frac{(2n-1)!!}{2^{n}} \sqrt{\pi}
\label{gamma8}
\end{equation}
\noindent
and $\mathbf{3.461.3}$ to 
\begin{equation}
\ift t^{n}e^{-t} \, dt = n!. 
\end{equation}
\noindent
These are now evaluated by identifying them with $\Gamma(n + \tfrac{1}{2})$
and $\Gamma(n+1)$, respectively.  \\

A second way to introduce fake parameters is to shift the integral 
(\ref{gamma1}) via $s = t+b$ to produce 
\begin{equation}
\int_{b}^{\infty} (s-b)^{a-1}e^{- s \mu} \, ds = \mu^{-a} e^{-\mu b} 
\Gamma(a).
\end{equation}

\noindent
This appears as $\mathbf{3.382.2}$ in \cite{gr}. \\

\noindent
There are many more integrals in \cite{gr} that can be reduced to the 
gamma function. These will be reported in a future publication.

\bigskip

\noindent
{\bf Acknowledgments}. The author wishes to thank Luis Medina for a careful 
reading of an earlier version of the paper. The partial support of
$\text{NSF-DMS } 0409968$ is also acknowledged.

\end{document}